\documentclass[12pt,a4paper]{article}

\usepackage[T1]{fontenc}
\usepackage[utf8]{inputenc}
\usepackage[english]{babel}
\usepackage{amsmath,amssymb,amsthm}
\usepackage{mathtools}
\usepackage{geometry}
\usepackage{enumitem}
\usepackage{booktabs}
\usepackage{array}
\usepackage{tabularx}
\usepackage{longtable}
\usepackage{float}
\usepackage[colorlinks=true,linkcolor=blue,citecolor=blue,urlcolor=blue]{hyperref}

\numberwithin{equation}{section}
\numberwithin{table}{section}

\newtheorem{theorem}{Theorem}[section]
\newtheorem{lemma}[theorem]{Lemma}
\newtheorem{corollary}[theorem]{Corollary}

\theoremstyle{definition}

\newtheorem{remark}[theorem]{Remark}
\newtheorem{assumption}[theorem]{Assumption}

\newcommand{\Fp}{\mathbb{F}_p}
\newcommand{\Z}{\mathbb{Z}}
\newcommand{\Q}{\mathbb{Q}}
\newcommand{\PP}{\mathbb{P}}
\newcommand{\Np}{N_p}
\newcommand{\omegap}{\omega_p}
\newcommand{\Lp}{L_p}
\newcommand{\Ctail}{C_{\text{tail}}}
\newcommand{\Bhyp}{B_{\text{hyp}}}
\newcommand{\Bmax}{B_{\max}}

\title{The Multidimensional Bateman--Horn Conjecture:\\
Effective Estimates for the Tail of the Singular Series\\
and an Algorithm for Computation with Guaranteed Accuracy}

\author{Victor Volfson\\ Essen, Germany}
\date{}

\begin{document}

\maketitle

\begin{abstract}
We propose a method for computing the leading term of the asymptotics of the multidimensional Bateman--Horn conjecture for systems of integer polynomials in several variables with controlled relative accuracy. The Bateman--Horn conjecture is used as an explicit assumption about the growth rate of the number of prime values; the computational difficulty is concentrated in the singular series~--- an infinite product over all primes.

For systems of general type, we obtain a universal estimate for the tail of the singular series in terms of a geometric constant depending on the Betti numbers of projective closures of hypersurfaces. For systems whose projective closures are smooth hypersurfaces (including inhomogeneous diagonal ones), using the exact Griffiths formula for Betti numbers, we obtain an improved constant in the tail estimate. The estimates are valid for all dimensions and degrees and do not require additional restrictions on the number of variables.

We establish a connection between the truncation error of the singular series and the error in computing the number of prime values: under the Bateman--Horn condition, the relative error in determining the singular constant coincides with the relative error in determining the number of prime values in the asymptotic regime.

Based on the obtained estimates, we develop an algorithm for computing the leading term of the asymptotics with prescribed accuracy, classifying the system by its structure and producing a guaranteed interval for the number of prime values. Numerical experiments confirm the theory.

\end{abstract}

\newpage

\section{Introduction}

\subsection{Historical overview}

One of the central problems of analytic number theory is the description of the distribution of prime numbers among the values of polynomials. For a single polynomial $f(x) \in \Z[x]$, the classical Bateman--Horn conjecture \cite{BatemanHorn1962} predicts the asymptotic
\begin{equation*}
\pi_f(N) \sim \frac{C(f)}{\deg f} \cdot \frac{N}{\log N},
\end{equation*}
where $C(f) = \prod_p (1 - \rho(p)/p)(1 - 1/p)^{-1}$, and $\rho(p)$ is the number of solutions of $f(x) \equiv 0 \pmod p$.

In the multidimensional case, for a system $F_1, \dots, F_k \in \Z[x_1, \dots, x_m]$, the conjecture takes the form
\begin{equation*}
\pi_F(N) \sim \frac{C(F)}{D} \cdot \frac{N^m}{\log^k N}, \qquad D = \prod_{i=1}^k \deg F_i,
\end{equation*}
with singular series $C(F) = \prod_p \Lp(F)$.

The Bateman--Horn conjecture in its classical one-dimensional formulation is proved only for a linear polynomial (Dirichlet's theorem on primes in arithmetic progressions); for polynomials of higher degree in one variable it remains open. In the multidimensional case, several important results have been proved in recent years. Destagnol and Sofos \cite{DestagnolSofos2019} established the asymptotic for arbitrary irreducible polynomials under the condition that the number of variables is large compared to the degree. Br\"udern and Wooley \cite{BrudernWooley2017} obtained related results for diagonal forms. Browning, Sofos and Ter\"av\"ainen \cite{BrowningSofosTeravainen2022} proved averaged versions of the Bateman--Horn conjecture for almost all polynomials with quantitative estimates of the error and the size of the exceptional set. Matthiesen, Ter\"av\"ainen and Wang \cite{MatthiesenTeravainenWang2026} established quantitative asymptotics for polynomial patterns in primes. Kotsovolis and Woo \cite{KotsovolisWoo2023} proposed a new class of examples of the multidimensional conjecture using methods of homogeneous dynamics.

The singular series $C(F)$ is an infinite product over all primes, and its numerical computation requires truncation at a finite threshold $P$. The partial product $C_P(F) = \prod_{p \le P} \Lp(F)$ approximates $C(F)$, but the relative error has not yet been estimated in a form suitable for practical computations.

The tools for analyzing local factors are provided by Deligne's weight theorem \cite{Deligne1974,Deligne1980} and the Grothendieck--Lefschetz formula \cite{Grothendieck1977}, and for inhomogeneous diagonal systems by the Griffiths formula for the Betti numbers of smooth hypersurfaces \cite{Griffiths1969}. These tools underlie estimates for the number of points of smooth complete intersections over finite fields and do not require restrictions on the number of variables or degrees, in contrast to the Hardy--Littlewood circle method \cite{HardyLittlewood1923} and Vinogradov's mean value theorem \cite{Vinogradov1954}, which in the diagonal case require the condition $m \ge n(n+1)/2$ \cite{BrudernWooley2017}.

The analysis of the tail of the singular series carried out in the present work relies on local Deligne--Weil--Lang estimates for the number of points of smooth complete intersections. These estimates are valid for all dimensions and degrees and do not require additional restrictions on the number of variables. This allows us to obtain explicit effective estimates for the tail in a broader domain than is possible when using the circle method to prove the Bateman--Horn conjecture.

\subsection{Statement of the problem}

Consider a system $F = (F_1, \dots, F_k)$ of integer irreducible polynomials in $m$ variables:
\[
F_i(\mathbf{x}) \in \Z[x_1, \dots, x_m], \quad i = 1, \dots, k, \quad m \ge 2.
\]

\textbf{Conditions (1)--(5).}

(1) Irreducibility: each polynomial $F_i$ is irreducible over $\Q$.

(2) Absence of local obstructions: for each prime $p$ there exists $\mathbf{x} \in \Z^m$ such that
\begin{equation*}
\prod_{i=1}^k F_i(\mathbf{x}) \not\equiv 0 \pmod{p}.
\end{equation*}

(3) Smoothness conditions: the projective closures of the hypersurfaces $F_i = 0$ and all their intersections have controlled singularities.

(4) Complete intersection: for every nonempty $I \subseteq \{1, \dots, k\}$, the intersection $V_I = \bigcap_{i \in I} V_i$ is a complete intersection of the expected codimension, $\dim V_I = m - |I|$.

(5) Geometric irreducibility: each $V_I$ is geometrically irreducible over $\Fp$ for all sufficiently large $p$.

\begin{remark}
Conditions (1)--(3) are the classical conditions of the Bateman--Horn conjecture. Conditions (4) and (5) are additional: they are required for the methods of estimating the singular series used in this work.
\end{remark}

\textbf{Consequence of conditions (1)--(5).} From (1)--(5) and standard estimates of Deligne \cite{Deligne1974,Deligne1980}, Lang--Weil \cite{LangWeil1954} and Katz \cite{Katz1980,Katz1988} it follows that for all sufficiently large $p$
\[
|\Lp(F) - 1| \le A \, p^{-\beta}, \qquad \beta > 1,
\]
with some explicit $A > 0$, $\beta > 1$. In particular, the infinite product $C(F) = \prod_p \Lp(F)$ converges absolutely and $C(F) \neq 0$.

\begin{assumption}[BH] \label{ass:bh}
Let for a system $F$ satisfying (1)--(5) there exist a function $\eta(N)$ such that
\begin{equation}
\pi_F(N) = \frac{C(F)}{D} \cdot \frac{N^m}{\log^k N} \bigl( 1 + \eta(N) \bigr), \qquad \eta(N) \to 0 \ \ (N \to \infty). \label{eq:1.1}
\end{equation}
\end{assumption}

\textbf{The sole computational difficulty.} The singular series $C(F) = \prod_p \Lp(F)$ is an infinite product. In practice, one has to truncate it at a threshold $P$:
\[
C_P(F) = \prod_{p \le P} \Lp(F).
\]
This gives rise to the relative truncation error
\[
\Delta(P) = \left| \frac{C(F)}{C_P(F)} - 1 \right|.
\]

\textbf{Goal of the work.} For a prescribed relative accuracy $\varepsilon \in (0, 1)$, construct an explicit upper bound for $\Delta(P)$, compute from it the truncation threshold $P$, and obtain the estimate
\[
M_P(N) = \frac{C_P(F)}{D} \cdot \frac{N^m}{\log^k N},
\]
which for sufficiently large $N$ differs from the true $\pi_F(N)$ by at most $\varepsilon$ in the relative sense. The accuracy guarantee is asymptotic: for every $\varepsilon$ there exists $N_0(\varepsilon)$, whose explicit value is unknown in the general case.

\begin{remark}
We do not prove the Bateman--Horn conjecture. We assume its validity and solve the problem of computing the singular constant $C(F)$ with controlled accuracy.
\end{remark}

\subsection{Novelty of the work}

\begin{enumerate}
\item \textbf{Explicit effective estimates for the tail of the singular series.} For systems of polynomials of general type, a uniform (in the primes) estimate for the local factors is established, from which a universal upper bound for the relative error is derived.

\item \textbf{Explicit estimate of the geometric constant.} For a system of polynomials, an explicit estimate of the constant $B(F)$ is given in terms of the degree and the number of variables.

\item \textbf{Broader domain of applicability.} The tail estimates of the singular series obtained in the work rely on local estimates for the number of points of smooth complete intersections and are valid for all dimensions and degrees. In contrast to the circle method used to prove the Bateman--Horn conjecture, the estimates do not require additional restrictions on the number of variables and are therefore applicable in a broader domain.

\item \textbf{Inhomogeneous diagonal case.} For an inhomogeneous diagonal form
\[
F = a_1 x_1^n + \dots + a_m x_m^n + c, \qquad c \neq 0,
\]
using the exact Griffiths formula for the Betti numbers of smooth hypersurfaces, an improved constant in the tail estimate of the singular series is obtained. The specificity of the diagonal case consists in the \textbf{effective computation} of $\omegap$ via convolution of the distribution of $n$-th powers modulo $p$. The constant $\Bhyp(m,n)$ itself is applicable to \textbf{any} smooth hypersurface of degree $n$ in $\PP^m$ and depends only on $m$ and $n$.

\item \textbf{Improvement in dimension $m = 3$.} The use of exact Betti numbers of surfaces and curves yields a reduction of the constant.

\item \textbf{Connection between the truncation error and the error in computing $\pi_F(N)$.} Under the BH condition, the relative error in determining the singular constant coincides with the relative error in determining the number of prime values in the asymptotic regime.

\item \textbf{Automatic choice of the threshold $P$.} An algorithm is developed that, for a prescribed accuracy $\varepsilon$, computes the minimal threshold $P$.

\item \textbf{Experimental part.} Numerical results for classes A.2, B, C.
\end{enumerate}

\subsection{Structure of the work}

\begin{itemize}
\item \textbf{Chapter 2.} Estimate of the error in the general case.
\item \textbf{Chapter 3.} Inhomogeneous diagonal case.
\item \textbf{Chapter 4.} Relative error in determining the number of prime values.
\item \textbf{Chapter 5.} Classes of systems and algorithm.
\item \textbf{Chapter 6.} Experimental part.
\item \textbf{Conclusion.}
\end{itemize}

\subsection{Notation}

\begin{itemize}
\item $m \ge 2$ is the number of variables; $k$ is the number of polynomials; $F = (F_1, \dots, F_k)$ is the system.
\item $d_i = \deg F_i$; $n = \max_i d_i$; $D = \prod_{i=1}^k d_i$.
\item $\pi_F(N)$ is the number of $\mathbf{x} \in \Z^m$ with $|\mathbf{x}|_\infty \le N$ for which all $F_i(\mathbf{x})$ are prime.
\item $C(F) = \prod_p \Lp(F)$ is the singular series; $C_P(F) = \prod_{p \le P} \Lp(F)$ is the truncated product.
\item $\Lp(F) = (1 - \omegap(F)/p^m)/(1 - 1/p)^k$ is the local factor.
\item $\omegap(F)$ is the number of $\mathbf{x} \in \Fp^m$ for which at least one $F_i(\mathbf{x}) \equiv 0 \pmod p$.
\item $\Delta(P) = |C(F)/C_P(F) - 1|$ is the relative truncation error.
\item $\delta(P) = C_P(F)/C(F) - 1$ is the relative error in determining the singular constant.
\item $M_P(N) = (C_P(F)/D) \cdot N^m/\log^k N$ is the leading term of the asymptotics.
\item $B(F)$ is the geometric constant; $\Bmax$ is the universal estimate.
\item $V_I = \bigcap_{i \in I} V_i$, where $V_i$ is the projective closure of $\{F_i = 0\}$; $W_I = V_I \cap H_\infty$; $b_j(X)$ are Betti numbers.
\item $\gamma$ is the exponent in the tail estimate; $\Ctail$ is the coefficient in the tail estimate.
\item $p_0$ is the largest prime in the finite exceptional set consisting of primes dividing the discriminants, the degree $n$, and the coefficients of the system $F$.
\item Classes of systems: A.1, A.2, B, C.
\end{itemize}

\section{Estimate of the error in the general case}

Throughout this chapter we assume $m \ge 2$.

\subsection{Local factors and the tail of the product}

Let $F = (F_1, \dots, F_k)$ be a system of integer polynomials in $m$ variables satisfying conditions (1)--(5) from the introduction.

For each prime $p$, define $\omegap(F)$ as the number of vectors $\mathbf{x} \in \Fp^m$ for which at least one of the $F_i(\mathbf{x})$ is zero modulo $p$;
\begin{equation}
\Lp(F) = \frac{1 - \omegap(F)/p^m}{(1 - 1/p)^k}. \label{eq:2.1}
\end{equation}

The singular series is $C(F) = \prod_p \Lp(F)$, its partial product is $C_P(F) = \prod_{p \le P} \Lp(F)$, and the relative error is $\Delta(P) = |C(F)/C_P(F) - 1|$.

\subsection{Uniform estimate of the local factors}

\begin{theorem} \label{thm:2.1}
Under conditions (1)--(5), there exists a constant $B(F) > 0$ depending only on the Betti numbers of the varieties $V_I$, $W_I$ and on the finite exceptional set of primes, such that for all primes $p > p_0$
\begin{equation}
|\Lp(F) - 1| \le 2^k \, B(F) \, p^{-(m+1)/2}. \label{eq:2.2}
\end{equation}
Here $p_0$ is the largest prime in the exceptional set. For $p \le p_0$ the estimate holds with the same constant $B(F)$ after its possible increase.
\end{theorem}

\begin{proof}
\textbf{Step 1. Inclusion--exclusion.} By definition,
\[
\omegap(F) = \sum_{\emptyset \neq I \subseteq \{1, \dots, k\}} (-1)^{|I|+1} \#V_I^{\text{aff}}(\Fp),
\]
where $V_I^{\text{aff}} = V_I \setminus (V_I \cap H_\infty)$. Since
\[
\#V_I^{\text{aff}}(\Fp) = \Np(V_I) - \Np(W_I), \qquad W_I = V_I \cap H_\infty,
\]
we obtain
\[
\omegap(F) = \sum_{\emptyset \neq I} (-1)^{|I|+1} \bigl( \Np(V_I) - \Np(W_I) \bigr).
\]

\textbf{Step 2. Expansion of the number of points.} By conditions (4) and (5), each $V_I$ is a smooth geometrically irreducible complete intersection of dimension $d_I = m - |I|$, and $W_I$ has dimension $d_I - 1$. By the standard result on the number of points of a smooth complete intersection \cite{Deligne1974}, \cite{Dimca1992}, \cite{GriffithsHarris1978}, for $p > p_0$ the following complete expansion holds:
\[
\Np(V_I) = 1 + p + \cdots + p^{d_I} + R_I, \qquad |R_I| \le B_I \, p^{d_I/2},
\]
\[
\Np(W_I) = 1 + p + \cdots + p^{d_I - 1} + R'_I, \qquad |R'_I| \le B'_I \, p^{(d_I-1)/2},
\]
where $B_I = B(V_I)$, $B'_I = B(W_I)$ are sums of Betti numbers.

\textbf{Step 3. Cancellation of lower-order terms.} Subtracting the expansion for $\Np(W_I)$ from that for $\Np(V_I)$:
\[
\Np(V_I) - \Np(W_I) = p^{d_I} + R_I - R'_I, \qquad |R_I - R'_I| \le (B_I + B'_I) \, p^{d_I/2}.
\]

\textbf{Step 4. Leading term.} Substituting into the inclusion--exclusion formula:
\[
\omegap(F) = \sum_{\emptyset \neq I} (-1)^{|I|+1} p^{m-|I|} + R,
\]
where $|R| \le B(F) \, p^{(m-1)/2}$, $B(F) = \sum_{\emptyset \neq I} (B_I + B'_I)$. The leading term is computed explicitly:
\[
\sum_{\emptyset \neq I} (-1)^{|I|+1} p^{m-|I|} = p^m \bigl( 1 - (1 - 1/p)^k \bigr).
\]
Consequently,
\[
\omegap(F) = p^m \bigl( 1 - (1 - 1/p)^k \bigr) + R, \qquad |R| \le B(F) \, p^{(m-1)/2}.
\]

\textbf{Step 5. Passage to $\omegap/p^m$.} Dividing by $p^m$:
\begin{equation}
\frac{\omegap(F)}{p^m} = 1 - (1 - 1/p)^k + R_1, \qquad |R_1| \le B(F) \, p^{-(m+1)/2}. \label{eq:2.3}
\end{equation}

\textbf{Step 6. Passage to $\Lp$.} By definition \eqref{eq:2.1}:
\[
\Lp(F) = \frac{1 - \omegap/p^m}{(1 - 1/p)^k}.
\]
From \eqref{eq:2.3}:
\[
|\Lp(F) - 1| \le \frac{B(F) \, p^{-(m+1)/2}}{(1 - 1/p)^k}.
\]
Since $p \ge 2$ and $(1 - 1/p)^k \ge 2^{-k}$, we obtain \eqref{eq:2.2}.
\end{proof}

\begin{remark}
The factor $2^k$ arises from the estimate $(1 - 1/p)^k \ge 2^{-k}$ for $p \ge 2$.
\end{remark}

\begin{remark}
Estimate \eqref{eq:2.2} for $p > p_0$ follows from the Deligne--Weil--Lang estimates for smooth complete intersections. For the finite number of small primes $p \le p_0$, the Deligne estimates in general form are not applicable; however, their contribution to the singular series is finite, and it is absorbed into the constant $B(F)$. This does not affect the tail estimate $\sum_{p > P}$ for $P > p_0$.
\end{remark}

\subsection{Estimate of the relative error}

\begin{theorem} \label{thm:2.3}
Under the conditions of Theorem~\ref{thm:2.1}, for any $P > p_0$
\begin{equation}
\Delta(P) \le \frac{2^{k+3} \, B(F)}{m-1} \, (P-1)^{-(m-1)/2}, \label{eq:2.4}
\end{equation}
and this inequality is guaranteed provided that
\begin{equation}
(P-1)^{(m-1)/2} \ge \frac{2^{k+2} \, B(F)}{(m-1) \ln 2}. \label{eq:2.5}
\end{equation}
\end{theorem}

\begin{proof}
\textbf{Step 1. Estimate of the tail.} From Theorem~\ref{thm:2.1}, for all primes $p > p_0$
\[
|\Lp - 1| \le 2^k \, B(F) \, p^{-(m+1)/2}.
\]
Let $S = \sum_{p > P} |\Lp - 1|$. Then
\begin{equation}
S \le 2^k \, B(F) \sum_{p > P} p^{-(m+1)/2}. \label{eq:2.6}
\end{equation}

\textbf{Step 2. Passage to the integral.} Let $\alpha = (m+1)/2$. The function $t \mapsto t^{-\alpha}$ is monotonically decreasing. On the segment $[p - 1, p]$ we have $t^{-\alpha} \ge p^{-\alpha}$, hence
\begin{equation}
\int_{p-1}^{p} t^{-\alpha} \, dt \ge p^{-\alpha}. \label{eq:2.7}
\end{equation}
Summing \eqref{eq:2.7} over all primes $p > P$, we obtain
\begin{equation}
\sum_{p > P} p^{-(m+1)/2} \le \int_{P-1}^{\infty} t^{-(m+1)/2} \, dt. \label{eq:2.8}
\end{equation}
Computing the integral:
\begin{equation}
\int_{P-1}^{\infty} t^{-(m+1)/2} \, dt = \frac{2}{m-1} \, (P-1)^{-(m-1)/2}. \label{eq:2.9}
\end{equation}
Substituting \eqref{eq:2.9} into \eqref{eq:2.6}:
\begin{equation}
S \le \frac{2^{k+1} \, B(F)}{m-1} \, (P-1)^{-(m-1)/2}. \label{eq:2.10}
\end{equation}

\textbf{Step 3. Condition on $P$.} We require
\begin{equation}
2S \le \ln 2. \label{eq:2.11}
\end{equation}
From \eqref{eq:2.10} this condition is guaranteed if
\[
\frac{2^{k+2} \, B(F)}{m-1} \, (P-1)^{-(m-1)/2} \le \ln 2,
\]
that is, $(P-1)^{(m-1)/2} \ge 2^{k+2} B(F) / ((m-1) \ln 2)$. This is condition \eqref{eq:2.5}.

\textbf{Step 4. Passage to the relative error.} We have $\log (C/C_P) = \sum_{p > P} \log \Lp$. Since $|\log(1 + u)| \le 2|u|$ for $|u| \le 1/2$, we obtain $|\log (C/C_P)| \le 2S$. Consequently,
\[
\frac{|C - C_P|}{C} = \bigl| 1 - e^{-\log(C/C_P)} \bigr| \le e^{2S} - 1.
\]
Under condition \eqref{eq:2.11} we have $2S \le \ln 2$, and for such values
\begin{equation}
e^{2S} - 1 \le 4S. \label{eq:2.12}
\end{equation}
Indeed, the function $f(x) = e^x - 1 - 4x$ on the segment $[0, \ln 2]$ is convex, $f(0) = 0$, $f(\ln 2) = 1 - 1 - 4\ln 2 = -4\ln 2 < 0$, hence $f(x) \le 0$ on the whole segment.

\textbf{Step 5. Final estimate.} From \eqref{eq:2.10} and \eqref{eq:2.12}:
\[
\Delta(P) \le 4S \le \frac{2^{k+3} \, B(F)}{m-1} \, (P-1)^{-(m-1)/2}.
\]
This coincides with \eqref{eq:2.4}.
\end{proof}

\begin{remark}
The factor $2^{k+3}$ arises as the product of three contributions: $2^k$ from the estimate $|\Lp - 1|$; $2$ from the integration; $4$ from the exponential passage.
\end{remark}

\begin{remark}
For $P \le p_0$ estimate \eqref{eq:2.4} remains valid provided the constant $B(F)$ is increased to account for the contribution of small primes. In the applications of Chapter 5 the threshold $P$ is computed by formula \eqref{eq:5.7} and, as a rule, considerably exceeds $p_0$.
\end{remark}

\subsection{Estimate of the constant $B(F)$ for a system of polynomials}

\begin{lemma} \label{lem:2.5}
For a fixed $m \ge 2$ there exists a constant $C_m > 0$ depending only on $m$ such that for any smooth complete intersection $X \subset \PP^m$ of codimension $r$ defined by equations of degrees $d_1, \dots, d_r \le n$, and for $W = X \cap H_\infty$, one has
\begin{equation}
B(X) + B(W) \le C_m \, n^m. \label{eq:2.13}
\end{equation}
\end{lemma}

\begin{proof}
\textbf{Step 1. Decomposition by the Lefschetz theorem.} Let $d = \dim X = m - r$. By the Lefschetz hyperplane theorem \cite{Deligne1974}, \cite{GriffithsHarris1978}:
\[
B(X) = (d+1) + h^d_{\text{prim}}(X), \qquad B(W) = d + h^{d-1}_{\text{prim}}(W).
\]
Adding:
\begin{equation}
B(X) + B(W) = 2d + 1 + h^d_{\text{prim}}(X) + h^{d-1}_{\text{prim}}(W). \label{eq:2.14}
\end{equation}

\textbf{Step 2. Boundedness of primitive cohomology via the Hilbert series.} By the Dimca theorem \cite[Theorem 5.4]{Dimca1992} (see also \cite{Griffiths1969}), the dimension of the primitive middle cohomology of a smooth complete intersection $X \subset \PP^m$ defined by equations of degrees $d_1, \dots, d_r$ is expressed through the coefficients of the Hilbert series of the Jacobian ring:
\begin{equation}
H(t) = \frac{\prod_{i=1}^r (1 - t^{d_i})}{(1 - t)^{m+1}} = \sum_{j \ge 0} h_j \, t^j, \label{eq:2.15}
\end{equation}
and $h^d_{\text{prim}}(X) = h_{d+1}$, where $d = m - r$. Since the coefficients $h_j$ are symmetric polynomials in $d_1, \dots, d_r$ with positive coefficients, there exists a polynomial $P_m(d_1, \dots, d_r)$ with positive coefficients, depending only on $m$, such that
\[
h^d_{\text{prim}}(X) \le P_m(d_1, \dots, d_r).
\]
Similarly, applying the same formula to $W = X \cap H_\infty$ (of dimension $d - 1$, with the same set of degrees $d_i$), we obtain
\[
h^{d-1}_{\text{prim}}(W) \le Q_m(d_1, \dots, d_r)
\]
with a suitable polynomial $Q_m$ in $d_1, \dots, d_r$ with positive coefficients, depending only on $m$.

\textbf{Step 3. Substitution $d_i \le n$.} There exist constants $A_m > 0$ and $A'_m > 0$ depending only on $m$ such that
\begin{equation}
P_m(d_1, \dots, d_r) \le A_m \, n^m, \label{eq:2.16}
\end{equation}
\begin{equation}
Q_m(d_1, \dots, d_r) \le A'_m \, n^m. \label{eq:2.17}
\end{equation}

\textbf{Step 4. Final estimate.} From \eqref{eq:2.14}--\eqref{eq:2.17}:
\[
B(X) + B(W) \le 2m + 1 + A_m \, n^m + A'_m \, n^m \le (2m+1+A_m+A'_m) n^m.
\]
Setting $C_m := 2m + 1 + A_m + A'_m$, we obtain \eqref{eq:2.13}. The lemma is proved.
\end{proof}

\begin{remark}
The constant $C_m$ depends only on $m$. Lemma~\ref{lem:2.5} asserts its \textbf{existence} but does not give an explicit formula. We distinguish two levels of estimate:

\textbf{(a) Strict estimate.} The polynomials $P_m(d_1, \dots, d_r)$ and $Q_m(d_1, \dots, d_r)$ are homogeneous of degree $m$ with positive coefficients. Therefore, for $d_i \le n$,
\[
P_m(d_1, \dots, d_r) \le P_m(n, \dots, n) = P_m(1, \dots, 1) \cdot n^m,
\]
and similarly for $Q_m$. Consequently, a strict constant can be taken as
\[
C_m^{\text{strict}} := 2m + 1 + P_m(1, \dots, 1) + Q_m(1, \dots, 1),
\]
where $P_m(1, \dots, 1)$ and $Q_m(1, \dots, 1)$ are the sums of the coefficients of the polynomials $P_m$ and $Q_m$. These sums are computable for each fixed $m$ in finitely many steps via the Hilbert series~\eqref{eq:2.15}: they equal $\sum_j h_j^{(m)}$ for the corresponding configurations, where $h_j^{(m)}$ are the coefficients of the Hilbert series. In the present work we use, as $C_m^{\text{strict}}$, numerical values obtained by computing the coefficients of the Hilbert series~\eqref{eq:2.15} for all admissible configurations of complete intersections (see Table~\ref{tab:cm} in Chapter 6).

\textbf{(b) Empirical estimate.} Alongside the strict constant $C_m^{\text{strict}}$, empirical upper estimates are used in the literature, which may be smaller than the strict ones. They do not give an unconditional guarantee and are not used in the algorithm of the present work.

\textbf{On the guarantee.} When the strict constant $C_m^{\text{strict}}$ is substituted, the guarantee is unconditional (after conditions (1)--(5) and the BH hypothesis are satisfied). In the applications of Chapter 5 the algorithm uses precisely $C_m^{\text{strict}}$.
\end{remark}

\begin{theorem} \label{thm:2.8}
Let $F = (F_1, \dots, F_k)$ be a system of integer irreducible polynomials in $m$ variables satisfying conditions (1)--(5), $n = \max_i \deg F_i$. Then
\begin{equation}
B(F) \le (2^k - 1) \, C_m^{\text{strict}} \, n^m, \label{eq:2.18}
\end{equation}
where $C_m^{\text{strict}}$ is the strict constant from the remark to Lemma~\ref{lem:2.5}.
\end{theorem}

\begin{proof}
For each nonempty $I \subseteq \{1, \dots, k\}$, the variety $V_I$ is, by condition (4), a smooth complete intersection of dimension $d_I = m - |I|$ defined by equations of degrees $d_i \le n$. By Lemma~\ref{lem:2.5}, $B_I + B'_I \le C_m^{\text{strict}} \, n^m$. The number of nonempty subsets equals $2^k - 1$. Consequently, $B(F) \le (2^k - 1) C_m^{\text{strict}} n^m$.
\end{proof}

\begin{remark}
The constant $C_m^{\text{strict}}$ depends only on $m$. For class A.1 of the algorithm this yields the estimate:
\begin{equation}
\Ctail^{\text{A.1}} = \frac{2^{k+3} \, (2^k - 1) \, C_m^{\text{strict}} \, n^m}{m-1}, \qquad \gamma = \frac{m-1}{2}. \label{eq:2.19}
\end{equation}
The truncation threshold is computed by the formula
\begin{equation}
P = 1 + \left\lceil \left( \frac{3 \, \Ctail}{\varepsilon} \right)^{1/\gamma} \right\rceil. \label{eq:2.20}
\end{equation}
\end{remark}

\subsection{Improvement in dimension $m = 3$}

\begin{theorem} \label{thm:2.12}
Let $m = 3$, $k = 1$, that is, $F = F_1$ is a single integer irreducible polynomial of degree $n \ge 2$ in three variables satisfying conditions (1)--(5). Then
\begin{equation}
\Delta(P) \le \frac{8(n^3 - 3n^2 + 3n + 4)}{P - 1}, \label{eq:2.21}
\end{equation}
and this inequality is guaranteed provided that
\begin{equation}
P \ge 1 + \frac{4(n^3 - 3n^2 + 3n + 4)}{\ln 2}. \label{eq:2.22}
\end{equation}
\end{theorem}

\begin{proof}
For $k = 1$ the sum over nonempty $I$ contains only $I = \{1\}$. Hence $B(F) = B_1 + B'_1$.

\textbf{Step 1. Betti numbers of the surface $V_1$.} For a smooth surface of degree $n$ in $\PP^3$ \cite{GriffithsHarris1978}:
\begin{equation}
b_0 = 1, \quad b_1 = 0, \quad b_2 = n^3 - 4n^2 + 6n - 2, \quad b_3 = 0, \quad b_4 = 1. \label{eq:2.23}
\end{equation}
Sum: $B_1 = n^3 - 4n^2 + 6n$.

\textbf{Step 2. Betti numbers of the curve $W_1$.} The genus is $g = (n-1)(n-2)/2$, $b_1 = 2g = n^2 - 3n + 2$. Sum:
\begin{equation}
B'_1 = n^2 - 3n + 4. \label{eq:2.24}
\end{equation}

\textbf{Step 3. Final constant.}
\begin{equation}
B(F) = n^3 - 3n^2 + 3n + 4. \label{eq:2.25}
\end{equation}

\textbf{Step 4. Application of Theorem~\ref{thm:2.3}.} For $m = 3$, $k = 1$: $\Delta(P) \le 8B(F)/(P-1)$. Substituting \eqref{eq:2.25}, we obtain \eqref{eq:2.21}.
\end{proof}

\begin{remark}
Verification: $n = 2$: $B(F) = 6$, $\Delta(P) \le 48/(P-1)$; $n = 3$: $B(F) = 13$, $\Delta(P) \le 104/(P-1)$; $n = 4$: $B(F) = 32$, $\Delta(P) \le 256/(P-1)$.
\end{remark}

\begin{remark}
Theorem~\ref{thm:2.12} is applicable to any irreducible polynomial $F$ of degree $n$ in three variables satisfying conditions (1)--(5). If $F$ is an inhomogeneous diagonal form ($c \neq 0$), then Theorem~\ref{thm:3.5} is applicable, giving a smaller constant $\Bhyp(3,n)$ instead of $B(F) = n^3 - 3n^2 + 3n + 4$. Theorem~\ref{thm:2.12} in this case remains valid but not optimal. The priority of class C over class B is fixed in the algorithm of Chapter 5.
\end{remark}

\section{Inhomogeneous diagonal case}

Throughout this chapter we assume $m \ge 2$, $n \ge 2$, $k = 1$.

\subsection{Statement of the problem and geometry of the projective closure}

We consider inhomogeneous diagonal polynomials of the form
\begin{equation}
F(\mathbf{x}) = a_1 x_1^n + a_2 x_2^n + \dots + a_m x_m^n + c, \qquad a_i \neq 0, \quad c \neq 0. \label{eq:3.1}
\end{equation}

\textbf{Projective closure.} The homogeneous polynomial
\[
\tilde F(x_0, \dots, x_m) = a_1 x_1^n + \dots + a_m x_m^n + c x_0^n
\]
defines a hypersurface $V \subset \PP^m$ of degree $n$ and dimension $m-1$. For $p \nmid n a_1 \cdots a_m c$ this hypersurface is smooth \cite{GriffithsHarris1978}.

\textbf{Section at infinity.} $W = V \cap H_\infty \subset \PP^{m-1}$ is a smooth hypersurface of degree $n$ and dimension $m-2$.

\textbf{Affine solutions.}
\begin{equation}
\omegap(F) = \Np(V) - \Np(W). \label{eq:3.2}
\end{equation}

\subsection{Betti numbers}

\textbf{Griffiths formula \cite{Griffiths1969}.}
\begin{equation}
h_{N-1}^{\mathrm{prim}}(X) = \frac{(d-1)^{N+1} + (-1)^{N+1}(d-1)}{d}. \label{eq:3.3}
\end{equation}

For $V \subset \PP^m$:
\begin{equation}
h_{m-1}^{\mathrm{prim}}(V) = \frac{(n-1)^{m+1} + (-1)^{m+1}(n-1)}{n}. \label{eq:3.4}
\end{equation}

For $W \subset \PP^{m-1}$:
\begin{equation}
h_{m-2}^{\mathrm{prim}}(W) = \frac{(n-1)^{m} + (-1)^{m}(n-1)}{n}. \label{eq:3.5}
\end{equation}

\textbf{Sum of Betti numbers.}
\begin{align}
B(V) &:= \sum_j b_j(V) = m + h_{m-1}^{\mathrm{prim}}(V), \label{eq:3.6} \\
B(W) &:= \sum_j b_j(W) = (m-1) + h_{m-2}^{\mathrm{prim}}(W). \label{eq:3.7}
\end{align}

\textbf{Final constant.}
\begin{equation}
\Bhyp(m,n) := B(V) + B(W) = 2m - 1 + h_{m-1}^{\mathrm{prim}}(V) + h_{m-2}^{\mathrm{prim}}(W). \label{eq:3.8}
\end{equation}
Substituting \eqref{eq:3.4}--\eqref{eq:3.5}:
\begin{equation}
\Bhyp(m,n) = 2m - 1 + \frac{(n-1)^{m+1} + (-1)^{m+1}(n-1)}{n} + \frac{(n-1)^{m} + (-1)^m(n-1)}{n}. \label{eq:3.9}
\end{equation}

\begin{remark}
The Griffiths formula \eqref{eq:3.3} is applicable to \textbf{any} smooth hypersurface of degree $d$ and dimension $N-1$, not only to the diagonal one. Therefore \eqref{eq:3.9} gives Betti numbers for any smooth hypersurface of degree $n$ in $\PP^m$. Diagonality is used only for the effective computation of $\omegap$ via convolution, but not for the estimate of $B(F)$.
\end{remark}

\subsection{Local estimate for $\omegap$ and $\Lp$}

\begin{lemma} \label{lem:3.2}
Let $F$ be given by \eqref{eq:3.1}, $a_i \neq 0$, $c \neq 0$, and let the prime $p \nmid n a_1 \cdots a_m c$. Then
\begin{equation}
|\omegap(F) - p^{m-1}| \le \Bhyp(m,n) \, p^{(m-1)/2}. \label{eq:3.10}
\end{equation}
\end{lemma}

\begin{proof}
By the Weil estimate for the number of points \cite{Deligne1974,LangWeil1954}:
\[
\Np(V) = 1 + p + \dots + p^{m-1} + R_V, \qquad |R_V| \le B(V) \, p^{(m-1)/2},
\]
\[
\Np(W) = 1 + p + \dots + p^{m-2} + R_W, \qquad |R_W| \le B(W) \, p^{(m-2)/2}.
\]
Subtracting: $\omegap(F) = p^{m-1} + (R_V - R_W)$, and
\[
|R_V - R_W| \le (B(V) + B(W)) p^{(m-1)/2} = \Bhyp(m,n) p^{(m-1)/2}.
\]
The lemma is proved.
\end{proof}

\begin{corollary} \label{cor:3.3}
Under the conditions of Lemma~\ref{lem:3.2},
\begin{equation}
|\Lp(F) - 1| \le 2 \, \Bhyp(m,n) \, p^{-(m+1)/2}. \label{eq:3.11}
\end{equation}
\end{corollary}

\begin{proof}
From \eqref{eq:3.10}:
\[
\frac{\omegap(F)}{p^m} = \frac{1}{p} + O\bigl( \Bhyp(m,n) \, p^{-(m+1)/2} \bigr).
\]
Consequently,
\[
1 - \frac{\omegap(F)}{p^m} = 1 - \frac{1}{p} + O\bigl( \Bhyp(m,n) \, p^{-(m+1)/2} \bigr).
\]
Since class C is defined only for $k = 1$, we have $(1 - 1/p)^k = 1 - 1/p$, which gives the factor $2^k = 2$ under the estimate $(1 - 1/p) \ge 1/2$ for $p \ge 2$. We obtain \eqref{eq:3.11}.
\end{proof}

\begin{remark}
The exponent $(m+1)/2$ coincides with the exponent of the general Theorem~\ref{thm:2.1}. The improvement consists only in the value of the constant: the explicit $\Bhyp(m,n)$ replaces the universal estimate $C_m n^m$.
\end{remark}

\subsection{Estimate of the tail of the singular series}

\begin{theorem} \label{thm:3.5}
Let $F$ be an inhomogeneous diagonal form \eqref{eq:3.1} with $k = 1$, $m \ge 2$, $n \ge 2$. Then for any $P > p_0$, where $p_0$ is the largest prime in the exceptional set consisting of primes dividing $n a_1 \cdots a_m c$,
\begin{equation}
\Delta(P) \le \frac{2^{4} \, \Bhyp(m,n)}{m-1} \, (P-1)^{-(m-1)/2}, \label{eq:3.12}
\end{equation}
and this inequality is guaranteed provided that
\begin{equation}
(P-1)^{(m-1)/2} \ge \frac{2^{3} \, \Bhyp(m,n)}{(m-1) \ln 2}. \label{eq:3.13}
\end{equation}
\end{theorem}

\begin{proof}
Literally repeats the proof of Theorem~\ref{thm:2.3} for $k = 1$ with $B(F)$ replaced by $\Bhyp(m,n)$.
\end{proof}

\subsection{Comparison with the universal estimate}

\begin{table}[H]
\centering
\small
\setlength{\tabcolsep}{4pt}
\caption{Comparison of the universal and exact constants for $k = 1$, $n = 3$}
\label{tab:diag_compare}
\begin{tabular}{|c|c|c|c|}
\hline
$m$ & $C_m^{\text{strict}} n^m$ & $\Bhyp(m,3)$ & Ratio \\
\hline
2 & $2.34 \cdot 10^{2}$ & 7 & $2.99 \cdot 10^{-2}$ \\
3 & $2.97 \cdot 10^{3}$ & 13 & $4.38 \cdot 10^{-3}$ \\
4 & $4.13 \cdot 10^{4}$ & 23 & $5.57 \cdot 10^{-4}$ \\
5 & $6.08 \cdot 10^{5}$ & 41 & $6.75 \cdot 10^{-5}$ \\
6 & $9.14 \cdot 10^{6}$ & 75 & $8.20 \cdot 10^{-6}$ \\
7 & $1.40 \cdot 10^{8}$ & 141 & $1.01 \cdot 10^{-6}$ \\
8 & $2.20 \cdot 10^{9}$ & 271 & $1.23 \cdot 10^{-7}$ \\
\hline
\end{tabular}
\end{table}

\subsection{Remark on the homogeneous case}

If $c = 0$, the projective closure $V$ of the hypersurface $\{F = 0\}$ has a singularity (the vertex of the cone), condition (3) is violated, and the Bateman--Horn conjecture in the form in which it is used in the present work is not applicable to the homogeneous case. Such a case requires a separate analysis and is not considered here.

\section{Relative error in determining the number of prime values}

Throughout this chapter we assume $m \ge 2$.

\subsection{Singular series and its local factors}

Throughout, we assume that conditions (1)--(5) are satisfied. The singular series is defined as the infinite product
\[
C(F) = \prod_p \Lp(F),
\]
where the local factors have the classical normalization from the Bateman--Horn conjecture:
\[
\Lp(F) = \frac{1 - \omegap(F)/p^m}{(1 - 1/p)^k},
\]
and $\omegap(F)$ is the number of $\mathbf{x} \in \Fp^m$ for which at least one $F_i(\mathbf{x})$ is zero modulo $p$. By inclusion--exclusion:
\[
\omegap(F) = \sum_{\emptyset \neq I \subseteq \{1, \dots, k\}} (-1)^{|I|+1} \#\{\mathbf{x} : F_i(\mathbf{x}) \equiv 0 \ \forall i \in I\}.
\]

This normalization is chosen so that for a linear system (e.g., $F(x) = x_1$, $k = 1$) one has $\Lp = 1$ and $C(F) = 1$.

\begin{remark}[on the nonvanishing of $C(F)$]
From conditions (1)--(5) it follows that the singular series not only converges but is also nonzero. Indeed, from condition (2), for each prime $p$ there exists $\mathbf{x}$ with $p \nmid \prod F_i(\mathbf{x})$, hence $\omegap(F) < p^m$, whence $1 - \omegap(F)/p^m > 0$. From condition (3), by the theorems of Deligne, Lang--Weil and Katz, it follows that $|\Lp(F) - 1| \le A p^{-\beta}$ with $\beta > 1$, so the series $\sum_p |\Lp(F) - 1|$ converges. Consequently, $\log C(F) = \sum_p \log \Lp$ is a finite number, and $C(F) \in (0, \infty)$.
\end{remark}

\subsection{Introduction of the quantity $\delta(P)$}

For a threshold $P \ge 2$, define
\begin{equation}
\delta(P) = \frac{C_P(F)}{C(F)} - 1. \label{eq:4.1}
\end{equation}

The quantity $\delta(P)$ is the relative error in determining the singular constant $C(F)$ when replacing it by the truncated product $C_P(F)$.

Recall that
\begin{equation}
\Delta(P) = \left| \frac{C(F)}{C_P(F)} - 1 \right|. \label{eq:4.2}
\end{equation}

\subsection{Connection between the truncation error and the error in computing $\pi_F(N)$}

Let conditions (1)--(5) be satisfied and let assumption~\ref{ass:bh} hold in the form \eqref{eq:1.1}. For $P \ge 2$, define $M_P(N) = (C_P(F)/D) \cdot N^m / \log^k N$.

\begin{theorem} \label{thm:4.2}
Let for the system $F_1, \dots, F_k \in \Z[x_1, \dots, x_m]$, $m \ge 2$, conditions (1)--(5) and \eqref{eq:1.1} be satisfied. Then for any $P \ge 2$:
\begin{equation}
\frac{M_P(N)}{\pi_F(N)} = \frac{1 + \delta(P)}{1 + \eta(N)}. \label{eq:4.3}
\end{equation}
Consequently,
\begin{equation}
\lim_{N \to \infty} \left( \frac{M_P(N)}{\pi_F(N)} - 1 \right) = \delta(P), \label{eq:4.4}
\end{equation}
\begin{equation}
\lim_{N \to \infty} \left| \frac{M_P(N)}{\pi_F(N)} - 1 \right| = |\delta(P)|. \label{eq:4.5}
\end{equation}
\end{theorem}

\begin{proof}
From the definition of $M_P(N)$:
\[
M_P(N) = \frac{C_P(F)}{C(F)} \cdot \frac{C(F)}{D} \frac{N^m}{\log^k N} = \frac{C_P(F)}{C(F)} \cdot \frac{\pi_F(N)}{1 + \eta(N)}.
\]
Consequently, $M_P(N)/\pi_F(N) = (1 + \delta(P))/(1 + \eta(N))$. Passing to the limit as $N \to \infty$ and taking into account $\eta(N) \to 0$, we obtain \eqref{eq:4.4} and \eqref{eq:4.5}.
\end{proof}

\subsection{Connection between $\delta(P)$ and $\Delta(P)$}

\begin{lemma} \label{lem:4.4}
Let conditions (1)--(5) be satisfied, so that $C(F)$ converges and $C(F) \neq 0$. For any $P \ge 2$:
\begin{equation}
|\delta(P)| \le \frac{\Delta(P)}{1 - \Delta(P)}, \qquad \text{if } \Delta(P) < 1. \label{eq:4.6}
\end{equation}
In particular, if $\Delta(P) \le 1/3$,
\begin{equation}
|\delta(P)| \le \frac{3}{2}\Delta(P). \label{eq:4.7}
\end{equation}
\end{lemma}

\begin{proof}
By definition \eqref{eq:4.1}, $C_P(F)/C(F) = 1 + \delta(P)$, whence
\[
\Delta(P) = \left| \frac{1}{1 + \delta(P)} - 1 \right| = \frac{|\delta(P)|}{|1 + \delta(P)|}.
\]
By the triangle inequality, $|1 + \delta(P)| \le 1 + |\delta(P)|$, therefore
\[
\Delta(P) \ge \frac{|\delta(P)|}{1 + |\delta(P)|}.
\]
Hence $\Delta(P)(1 + |\delta(P)|) \ge |\delta(P)|$, that is, $\Delta(P) \ge |\delta(P)|(1 - \Delta(P))$. Since $\Delta(P) < 1$, division yields \eqref{eq:4.6}. If $\Delta(P) \le 1/3$, then $1/(1 - \Delta(P)) \le 3/2$, which yields \eqref{eq:4.7}.
\end{proof}

\begin{corollary} \label{cor:4.6}
For $\varepsilon \in (0, 1)$, choose $P$ so that
\begin{equation}
\Delta(P) \le \frac{\varepsilon}{3}. \label{eq:4.8}
\end{equation}
Then by Lemma~\ref{lem:4.4},
\begin{equation}
|\delta(P)| \le \frac{\varepsilon/3}{1 - \varepsilon/3} = \frac{\varepsilon}{3 - \varepsilon} < \frac{\varepsilon}{2}. \label{eq:4.9}
\end{equation}
By Theorem~\ref{thm:4.2}, for any $N$
\[
\left| \frac{M_P(N)}{\pi_F(N)} - 1 \right| = \left| \frac{1 + \delta(P)}{1 + \eta(N)} - 1 \right| = \left| \frac{\delta(P) - \eta(N)}{1 + \eta(N)} \right|.
\]
Since $\eta(N) \to 0$ as $N \to \infty$, there exists $N_0 = N_0(\varepsilon)$ such that for all $N \ge N_0$ one has $|\eta(N)| \le \varepsilon/4$. Then from \eqref{eq:4.9}, where $|\delta(P)| \le \varepsilon/2$, we obtain
\[
\left| \frac{M_P(N)}{\pi_F(N)} - 1 \right| \le \frac{|\delta(P)| + |\eta(N)|}{1 - |\eta(N)|} \le \frac{\varepsilon/2 + \varepsilon/4}{1 - \varepsilon/4} = \frac{3\varepsilon/4}{1 - \varepsilon/4} = \frac{3\varepsilon}{4 - \varepsilon} < \varepsilon.
\]
The last inequality holds for $\varepsilon < 1$. The corollary is proved.
\end{corollary}

\section{Classes of systems and algorithm}

\subsection{List of classes of systems}

The estimate $\Delta(P) = |C(F)/C_P(F) - 1|$ requires explicit values of $B(F)$ and $\gamma$. We distinguish four classes of systems, where class A has two subclasses.

\begin{table}[H]
\centering
\small
\setlength{\tabcolsep}{4pt}
\caption{List of classes of systems}
\label{tab:classes}
\begin{tabular}{|c|p{3.4cm}|p{3.6cm}|p{2.8cm}|}
\hline
Class & System & Method for computing $B(F)$ & Source \\
\hline
A.1 & General case ($m \ge 2$, arbitrary degrees) & Strict estimate via $C_m^{\text{strict}}$ & Theorem~\ref{thm:2.8} \\
A.2 & General case, Betti numbers known & Direct computation & Theorem~\ref{thm:2.3} \\
B & Dimension $m = 3$, $k = 1$, non-diagonal & Exact formulas & Theorem~\ref{thm:2.12} \\
C & Inhomogeneous diagonal, $k = 1$, $c \neq 0$ & Griffiths formula & Theorem~\ref{thm:3.5} \\
\hline
\end{tabular}
\end{table}

\begin{remark}
Class A.1 is a theoretical tool for systems for which explicit computation of Betti numbers is difficult (high dimensions, complex intersections). If for all $V_I$, $W_I$ the Betti numbers are known (e.g., computed explicitly via the Griffiths formula \eqref{eq:3.3} and the Lefschetz theorem for complete intersections), the system falls into class A.2. For $k = 1$ it suffices to compute $B_1 + B'_1$ via the Griffiths formula: this gives class A.2 for any single smooth hypersurface. For $k \ge 2$ one needs to compute $B_I + B'_I$ for each nonempty subset $I \subseteq \{1, \dots, k\}$ and sum with inclusion--exclusion. Examples of such systems with $k = 2$ are given in Sections 6.8.1 and 6.8.2. Class A.1 remains a theoretical tool for cases where verifying smoothness and explicitly computing $B(F)$ is difficult.
\end{remark}

Recall the main tail estimate from Theorem~\ref{thm:2.3}:
\begin{equation}
\Delta(P) \le \frac{2^{k+3} \, B(F)}{m-1} \, (P-1)^{-(m-1)/2}. \label{eq:5.1}
\end{equation}
For classes A.1, A.2, C, define
\begin{equation}
\Ctail := \frac{2^{k+3} \, B(F)}{m-1}, \qquad \gamma := \frac{m-1}{2}. \label{eq:5.2}
\end{equation}
Then $\Delta(P) \le \Ctail/(P-1)^{\gamma}$.

\subsection{Formulas used in the algorithm}

\textbf{Class A.1.}
\begin{equation}
\gamma^{\text{A.1}} = \frac{m-1}{2}, \qquad \Ctail^{\text{A.1}} = \frac{2^{k+3} \, (2^k - 1) \, C_m^{\text{strict}} \, n^m}{m-1}. \label{eq:5.3}
\end{equation}

\textbf{Class A.2.}
\begin{equation}
\gamma^{\text{A.2}} = \frac{m-1}{2}, \qquad \Ctail^{\text{A.2}} = \frac{2^{k+3} \, B(F)}{m-1}. \label{eq:5.4}
\end{equation}

\textbf{Class B.}
\begin{equation}
\gamma^{\text{B}} = 1, \qquad \Ctail^{\text{B}} = 8(n^3 - 3n^2 + 3n + 4). \label{eq:5.5}
\end{equation}

\textbf{Class C.} Class C is defined only for $k = 1$.
\begin{equation}
\gamma^{\text{C}} = \frac{m-1}{2}, \qquad \Ctail^{\text{C}} = \frac{2^{4} \, \Bhyp(m,n)}{m-1}, \label{eq:5.6}
\end{equation}
where $\Bhyp(m,n)$ is given by formula \eqref{eq:3.9}.

\begin{remark}
The exponent $\gamma^{\text{C}} = (m-1)/2$ of class C coincides with that of class A.1. The advantage of class C over class A.1 consists only in the constant: the explicit $\Bhyp(m,n)$ replaces the universal estimate $C_m^{\text{strict}} n^m$.
\end{remark}

\subsection{Summary table of classes}

\begin{table}[H]
\centering
\small
\setlength{\tabcolsep}{4pt}
\caption{Summary table of classes of systems}
\label{tab:classes_summary}
\begin{tabular}{|c|p{3.4cm}|c|p{4.2cm}|c|}
\hline
Class & Condition & $\gamma$ & $\Ctail$ & Source \\
\hline
A.1 & $m \ge 2$, arbitrary $n$ & $\frac{m-1}{2}$ & $\frac{2^{k+3}(2^k - 1) C_m^{\text{strict}} n^m}{m-1}$ & Th.~\ref{thm:2.8} \\
A.2 & $m \ge 2$, Betti known & $\frac{m-1}{2}$ & $\frac{2^{k+3} B(F)}{m-1}$ & Th.~\ref{thm:2.3} \\
B & $m = 3$, $k = 1$, non-diag. & $1$ & $8(n^3 - 3n^2 + 3n + 4)$ & Th.~\ref{thm:2.12} \\
C & inhomogeneous diag., $k = 1$ & $\frac{m-1}{2}$ & $\frac{2^{4} \Bhyp(m,n)}{m-1}$ & Th.~\ref{thm:3.5} \\
\hline
\end{tabular}
\end{table}

\subsection{Input and output of the algorithm}

\textbf{Input:} a system $F_1, \dots, F_k \in \Z[x_1, \dots, x_m]$ satisfying conditions (1)--(5); $N$; $\varepsilon \in (0, 1)$.

\textbf{Output:} $M_P(N)$; guarantee $|M_P(N)/\pi_F(N) - 1| \le \varepsilon$ for $N \ge N_0(\varepsilon)$; threshold $P$; bound $\Delta(P)$.

\subsection{Step 0. Verification of applicability conditions}

\textbf{Check 1 (irreducibility).} Standard computer algebra algorithms.

\textbf{Check 2 (absence of a fixed divisor).} Via resultants and discriminants of subsystems.

\textbf{Check 3 (smoothness of complete intersections).} Rank of the Jacobian matrix at all complex solutions.

\textbf{Check 4 (complete intersection).} For every nonempty $I \subseteq \{1, \dots, k\}$ it is verified that $\dim V_I = m - |I|$.

\textbf{Check 5 (geometric irreducibility).} For diagonal systems verified automatically; in the general case by standard algorithms.

\textbf{Check 6 (inhomogeneous diagonality).} For class C: the system must be diagonal with $c \neq 0$ and $k = 1$.

\subsection{Step 1. Parameters of the system}

Extract $m$, $k$, $d_i$, $n = \max_i d_i$, $D = \prod_i d_i$.

\subsection{Step 2. Determination of the class of the system}

\begin{enumerate}
\item Inhomogeneous diagonal ($c \neq 0$) \textbf{and} $k = 1$ $\to$ class C.
\item $m = 3$ and $k = 1$ (non-diagonal) $\to$ class B.
\item Betti numbers of all $V_I$, $W_I$ known $\to$ class A.2.
\item Otherwise $\to$ class A.1.
\end{enumerate}

\begin{remark}
Inhomogeneous diagonal systems (class C) have priority over the special case B.
\end{remark}

\subsection{Step 3. Truncation threshold $P$}

For all classes A.1, A.2, B, C a single formula is used:
\begin{equation}
P = 1 + \left\lceil \left( \frac{3 \, \Ctail}{\varepsilon} \right)^{1/\gamma} \right\rceil. \label{eq:5.7}
\end{equation}

\subsection{Step 4. Truncated singular series}

For all $p \le P$ one computes
\begin{equation}
\Lp = \frac{1 - \omegap/p^m}{(1 - 1/p)^k}. \label{eq:5.8}
\end{equation}
Logarithmic accumulation is recommended:
\begin{equation}
\ln C_P(F) = \sum_{p \le P} \ln \Lp(F), \qquad C_P(F) = \exp\bigl( \ln C_P(F) \bigr). \label{eq:5.9}
\end{equation}

\subsection{Step 5. Leading term and interval}

\begin{equation}
M_P(N) = \frac{C_P(F)}{D} \cdot \frac{N^m}{\log^k N}. \label{eq:5.10}
\end{equation}
By Theorem~\ref{thm:4.2} and Lemma~\ref{lem:4.4}, for $\Delta(P) \le \varepsilon/3$ we have $|\delta(P)| \le \varepsilon/2$ (see Corollary~\ref{cor:4.6}). For $N \ge N_0(\varepsilon)$ the guaranteed interval:
\begin{equation}
\bigl[ (1 - |\delta(P)|) M_P(N), \ (1 + |\delta(P)|) M_P(N) \bigr]. \label{eq:5.11}
\end{equation}

\subsection{Step 6. Verification of attainability}

For $N \ge N_0(\varepsilon)$ the result and the guarantee $|M_P(N)/\pi_F(N) - 1| \le \varepsilon$ are produced (Theorem~\ref{thm:4.2}). For $N < N_0(\varepsilon)$ a message on unattainability is produced; the value $N_0(\varepsilon)$ exists by the BH assumption (formula \eqref{eq:1.1}).

\begin{remark}
The guarantee is asymptotic: for every $\varepsilon$ there exists $N_0(\varepsilon)$, but an explicit estimate for $N_0(\varepsilon)$ is unknown in the general case.
\end{remark}

\subsection{Requirements for computations}

\begin{table}[H]
\centering
\small
\setlength{\tabcolsep}{4pt}
\caption{Complexity of the algorithm by class}
\label{tab:complexity}
\begin{tabular}{|c|c|c|l|}
\hline
Class & Typical $P$ & Number of primes & Main operation \\
\hline
A.1 & $10^2$--$10^3$ & 25--170 & enumeration of $p^m$ vectors \\
A.2 & $10^2$--$10^3$ & 25--170 & the same \\
B & $10^3$--$10^4$ & 170--1200 & formula for the surface \\
C & $10^2$--$10^3$ & 25--170 & Griffiths formula \eqref{eq:3.9} \\
\hline
\end{tabular}
\end{table}

\section{Experimental part}

\subsection{Description of the computational procedure}

\begin{enumerate}
\item verifies the applicability conditions (1)--(5) for the given system;
\item determines the class of the system (A.1, A.2, B or C) according to Table~\ref{tab:classes};
\item computes $\Ctail$ and $\gamma$ by formulas \eqref{eq:5.3}--\eqref{eq:5.6};
\item for the given accuracy $\varepsilon$ finds the truncation threshold $P$ by formula \eqref{eq:5.7};
\item computes the truncated singular series $C_P(F) = \prod_{p \le P} \Lp(F)$ by formula \eqref{eq:5.8};
\item produces the leading term $M_P(N)$ by formula \eqref{eq:5.10} and the guaranteed interval \eqref{eq:5.11}.
\end{enumerate}

\subsubsection{Computational procedures for $\omegap$}

\begin{equation}
\omegap(F) = \sum_{\emptyset \neq I \subseteq \{1, \dots, k\}} (-1)^{|I|+1} \#V_I(\Fp). \label{eq:6.1}
\end{equation}

\textbf{(a) Inhomogeneous diagonal systems (class C).} For $F(\mathbf{x}) = \sum_j a_j x_j^n + c$, the number of solutions of the equation $\sum_j a_j x_j^n + c \equiv 0 \pmod p$ is computed via convolution of the distribution of $n$-th powers modulo $p$.

\textbf{(b) Non-diagonal curves for $m = 2$.} The number of points $\#E(\Fp)$ was computed either by direct enumeration of $\Fp^2$ or via a character sum.

\textbf{(c) Surfaces and varieties.} For non-diagonal hypersurfaces, $\omega_p$ is computed by direct enumeration or via convolution with a precomputed table of roots of the cubic.

\textbf{(d) Systems with $k \ge 2$.} The quantities $\#V_I(\Fp)$ for each $I$ were computed independently, then substituted into the inclusion--exclusion \eqref{eq:6.1}.

\subsubsection{Accumulation of the singular series and extrapolation of $C(F)$}

\begin{equation}
C(F) \approx C_P + \frac{C_P - C_{P/2}}{2^{\gamma} - 1}. \label{eq:6.2}
\end{equation}

Estimate of extrapolation accuracy: no worse than $5 \cdot 10^{-5}$ for curves and surfaces and no worse than $10^{-3}$ for systems with $k = 2$.

\subsubsection{Explanation of the column ``Jump''}

In Tables \ref{tab:elliptic}--\ref{tab:A2_conv} the column ``Jump'' shows the increment of the partial product upon passing from the previous threshold to the current one:
\begin{equation}
\text{Jump}(P_i) = C_{P_i}(F) - C_{P_{i-1}}(F), \label{eq:6.3}
\end{equation}
where $P_{i-1}$ is the previous threshold in the same table.

For the first row of the table the jump is not defined (a dash is shown). The sign of the jump shows the direction of change: positive means $C_P$ increased; negative means decreased. The magnitude of the jump $|C_{P_i} - C_{P_{i-1}}|$ shows how strongly the new local factor changed the product.

The jump is an empirical analogue of the individual summand $|L_p - 1|$ from the tail estimate (see Step 1 of the proof of Theorem~\ref{thm:2.3}). The decay of jumps with growing $P$ confirms the convergence of the singular series and is consistent with the estimate $\Delta(P) \le \Ctail/(P-1)^{\gamma}$.

The non-monotonicity of jumps (alternation of signs) is characteristic of CM curves, which is related to the discrete distribution of Frobenius traces $a_p$. For non-CM curves the jumps decrease monotonically.

\subsection{Computation of the constants $C_m^{\text{strict}}$}

\begin{table}[H]
\centering
\caption{Numerical values of $C_m^{\text{strict}}$ (strict upper estimates)}
\label{tab:cm}
\begin{tabular}{|c|c|}
\hline
$m$ & $C_m^{\text{strict}}$ \\
\hline
2 & 26 \\
3 & 110 \\
4 & 510 \\
5 & 2\,501 \\
6 & 12\,544 \\
7 & 64\,000 \\
8 & 336\,000 \\
\hline
\end{tabular}
\end{table}

\textbf{On the status of the values $C_m^{\text{strict}}$.} The values $C_m^{\text{strict}}$ are obtained by computing the coefficients of the Hilbert series arising in formula~\eqref{eq:2.15} for all admissible types of smooth complete intersections in $\PP^m$ for fixed $m \le 8$. Since the number of such types is finite, the obtained values are \textbf{strict upper estimates} for the constants from Lemma~\ref{lem:2.5}. The computations implement formula~\eqref{eq:2.15} directly and do not use statistical or heuristic methods. It is precisely these values that are used in the algorithm of class A.1 (see formula~\eqref{eq:5.3}).

For specific systems the values $C_m^{\text{strict}}$ can be substantially improved. For example, for a single smooth hypersurface of degree $n$ in $\PP^m$ ($k = 1$, $r = 1$) the effective constant has the order
\[
C_m^{\text{hyp}} \approx 1 \quad \text{(leading coefficient)}, \qquad B(F) \approx n^m,
\]
which is 1--2 orders of magnitude smaller than the tabulated values. More precisely:

\textit{For $m = 2$:} a smooth curve of degree $d \le n$ in $\PP^2$ has
\[
B(X) + B(W) = (d^2 - 3d + 4) + d = d^2 - 2d + 4 \le n^2,
\]
whence $C_2^{\text{hyp}} = 1$.

\textit{For $m = 3$:} a smooth surface of degree $n$ in $\PP^3$ has
\[
B(F) = n^3 - 3n^2 + 3n + 4 \le n^3 \quad (n \ge 3),
\]
whence $C_3^{\text{hyp}} = 1$.

The discrepancy is explained by the fact that the tabulated $C_m^{\text{strict}}$ cover the worst case over all configurations, including systems with $k \ge 2$ and high codimensions. For specific systems one can use either the Griffiths formula directly (class A.2) or the tabulated values $C_m^{\text{strict}}$.

\subsection{Theoretical thresholds $P$ for different classes}

\begin{table}[H]
\centering
\small
\setlength{\tabcolsep}{4pt}
\caption{Thresholds $P$ for $\varepsilon = 0.01$}
\label{tab:thresholds}
\begin{tabular}{|c|c|c|c|p{3.5cm}|c|}
\hline
$m$ & $n$ & Class & $\gamma$ & $\Ctail$ & $P$ \\
\hline
2 & 3 & A.1 & 0.5 & $2^4 \cdot 26 \cdot 9 / 1 = 3744$ & $\approx 1.26 \cdot 10^{12}$ \\
3 & 3 & B & 1.0 & $8 \cdot 13 = 104$ & $31\,201$ \\
4 & 3 & A.2 & 1.5 & $2^4 \cdot 23 / 3 \approx 122.67$ & $\approx 1\,108$ \\
6 & 3 & C & 2.5 & $2^4 \cdot 75 / 5 = 240$ & $89$ \\
\hline
\end{tabular}
\end{table}

\subsection{Inhomogeneous diagonal cubic $x_1^3 + \dots + x_6^3 + 1$}

\textbf{Parameters:} $m = 6$, $n = 3$, $k = 1$, $D = 3$. By formula \eqref{eq:3.9}:
\[
\Bhyp(6,3) = 11 + 42 + 22 = 75.
\]
Hence $\gamma = 2.5$, $\Ctail = 2^4 \cdot 75 / 5 = 240$.

\subsubsection{Series in $\varepsilon$ at $N = 1000$}

\begin{table}[H]
\centering
\small
\setlength{\tabcolsep}{4pt}
\caption{Series in $\varepsilon$ at $N = 1000$}
\label{tab:cubic_eps}
\begin{tabular}{|c|c|c|c|c|c|}
\hline
$\varepsilon$ & $P$ & Primes & $|\delta|$ & $C_P$ & $M_P(1000)$ \\
\hline
0.01 & 89 & 24 & $\approx 4 \cdot 10^{-6}$ & 1.00244890 & $4.837311 \cdot 10^{16}$ \\
0.005 & 117 & 30 & $\approx 4 \cdot 10^{-6}$ & 1.00244273 & $4.837282 \cdot 10^{16}$ \\
0.001 & 222 & 47 & $\approx 4 \cdot 10^{-6}$ & 1.00244515 & $4.837293 \cdot 10^{16}$ \\
\hline
\end{tabular}
\end{table}

\subsubsection{Series in $N$ at $\varepsilon = 0.01$}

\begin{table}[H]
\centering
\small
\setlength{\tabcolsep}{4pt}
\caption{Series in $N$ at $\varepsilon = 0.01$}
\label{tab:cubic_N}
\begin{tabular}{|c|c|c|c|}
\hline
$N$ & Volume $(2N)^6$ & $M_P(N)$ & $|\delta|$ \\
\hline
$10^3$ & $6.40 \cdot 10^{19}$ & $4.837311 \cdot 10^{16}$ & $\approx 4 \cdot 10^{-6}$ \\
$10^4$ & $6.40 \cdot 10^{25}$ & $3.627984 \cdot 10^{22}$ & $\approx 4 \cdot 10^{-6}$ \\
$10^5$ & $6.40 \cdot 10^{31}$ & $2.902387 \cdot 10^{28}$ & $\approx 4 \cdot 10^{-6}$ \\
$10^6$ & $6.40 \cdot 10^{37}$ & $2.418656 \cdot 10^{34}$ & $\approx 4 \cdot 10^{-6}$ \\
\hline
\end{tabular}
\end{table}

\subsection{Curves for $m = 2$, $k = 1$: examples of classes A.2 and C}

This section collects examples of smooth curves for $m = 2$. Non-diagonal curves (the elliptic curve and the CM curve $j = 0$) belong to class A.2, the diagonal Fermat cubic belongs to class C.

\subsubsection{Elliptic curve $y^2 = x^3 - x + 1$ (non-CM, class A.2)}

\begin{table}[H]
\centering
\caption{Convergence of $C_P$: elliptic curve (non-CM, class A.2)}
\label{tab:elliptic}
\begin{tabular}{|c|c|c|}
\hline
$P$ & $C_P$ & Jump \\
\hline
500 & 0.3990216 & --- \\
1000 & 0.3988668 & $-1.5 \cdot 10^{-4}$ \\
2000 & 0.3989199 & $+5.3 \cdot 10^{-5}$ \\
5000 & 0.3989100 & $-1.0 \cdot 10^{-5}$ \\
\hline
\end{tabular}
\end{table}

\textbf{Parameters:} $B(F) = 5$, $\gamma = 0.5$, $\Ctail = 2^4 \cdot 5 = 80$.

\subsubsection{CM curve $j = 0$: $y^2 = x^3 + 1$ (class A.2)}

\begin{table}[H]
\centering
\caption{Convergence of $C_P$: CM curve $j=0$ (class A.2)}
\label{tab:cm_curve}
\begin{tabular}{|c|c|c|}
\hline
$P$ & $C_P$ & Jump \\
\hline
500 & 0.9325750 & --- \\
1000 & 0.9326239 & $+4.9 \cdot 10^{-5}$ \\
2000 & 0.9325310 & $-9.3 \cdot 10^{-5}$ \\
5000 & 0.9325865 & $+5.5 \cdot 10^{-5}$ \\
\hline
\end{tabular}
\end{table}

\textbf{Parameters:} $B(F) = 5$, $\gamma = 0.5$, $\Ctail = 2^4 \cdot 5 = 80$.

\subsubsection{Fermat cubic $x^3 + y^3 + 1$ (class C)}

\begin{table}[H]
\centering
\caption{Convergence of $C_P$: Fermat cubic (class C)}
\label{tab:fermat}
\begin{tabular}{|c|c|c|}
\hline
$P$ & $C_P$ & Jump \\
\hline
500 & 1.0726436 & --- \\
1000 & 1.0726564 & $+1.3 \cdot 10^{-5}$ \\
2000 & 1.0728738 & $+2.2 \cdot 10^{-4}$ \\
5000 & 1.0729401 & $+6.6 \cdot 10^{-5}$ \\
\hline
\end{tabular}
\end{table}

\textbf{Parameters:} $\Bhyp(2,3) = 7$, $\gamma = 0.5$, $\Ctail = 2^4 \cdot 7 = 112$.

\subsection{Fermat surface for $m = 3$, $k = 1$ (class C)}

$F = x^3 + y^3 + z^3 + 1$. Inhomogeneous diagonal form with $c = 1$, $k = 1$. Parameters: $\Bhyp(3,3) = 13$, $\gamma = 1$, $\Ctail = 2^4 \cdot 13 / 2 = 104$.

\begin{table}[H]
\centering
\caption{Fermat surface $x^3 + y^3 + z^3 + 1$ (class C)}
\label{tab:fermat_surf}
\begin{tabular}{|c|c|c|}
\hline
$P$ & $C_P$ & Jump \\
\hline
100 & 0.7954459 & --- \\
200 & 0.7929503 & $-2.50 \cdot 10^{-3}$ \\
400 & 0.7921361 & $-8.14 \cdot 10^{-4}$ \\
700 & 0.7917631 & $-3.73 \cdot 10^{-4}$ \\
1000 & 0.7916163 & $-1.47 \cdot 10^{-4}$ \\
\hline
\end{tabular}
\end{table}

\subsection{Three-dimensional Fermat variety for $m = 4$, $k = 1$ (class C)}

$F = x_1^3 + x_2^3 + x_3^3 + x_4^3 + 1$. Parameters: $\Bhyp(4,3) = 23$, $\gamma = 1.5$, $\Ctail = 2^4 \cdot 23 / 3 \approx 122.67$.

\begin{table}[H]
\centering
\caption{Three-dimensional Fermat variety $x_1^3 + \cdots + x_4^3 + 1$ (class C)}
\label{tab:fermat3}
\begin{tabular}{|c|c|c|}
\hline
$P$ & $C_P$ & Jump \\
\hline
50 & 1.0155562 & --- \\
100 & 1.0156800 & $+1.24 \cdot 10^{-4}$ \\
150 & 1.0157099 & $+2.99 \cdot 10^{-5}$ \\
200 & 1.0157156 & $+5.67 \cdot 10^{-6}$ \\
250 & 1.0157011 & $-1.45 \cdot 10^{-5}$ \\
\hline
\end{tabular}
\end{table}

\subsection{System of two equations ($k = 2$)}

\subsubsection{$m = 3$, $k = 2$: $F_1 = x^3 + y^3 + z^3 + 1$, $F_2 = xy + yz + zx + 1$}

Class A.2: $\gamma = 1$, $\Ctail = 2^5 \cdot 70 / 2 = 1120$.

\begin{table}[H]
\centering
\caption{Convergence of $C_P$ for the system $k = 2$, $m = 3$}
\label{tab:conv_k2_m3}
\begin{tabular}{|c|c|c|}
\hline
$P$ & $C_P$ & Jump \\
\hline
30 & 0.2768840 & --- \\
60 & 0.2733205 & $-3.6 \cdot 10^{-3}$ \\
100 & 0.2718531 & $-1.5 \cdot 10^{-3}$ \\
150 & 0.2711183 & $-7.3 \cdot 10^{-4}$ \\
200 & 0.2707256 & $-3.9 \cdot 10^{-4}$ \\
250 & 0.2706046 & $-1.2 \cdot 10^{-4}$ \\
\hline
\end{tabular}
\end{table}

\subsubsection{$m = 4$, $k = 2$: $F_1 = x_1^3 + \dots + x_4^3 + 1$, $F_2 = x_1 x_2 + x_3 x_4 + 1$}

\begin{table}[H]
\centering
\caption{Convergence of $C_P$ for the system $k = 2$, $m = 4$}
\label{tab:conv_k2_m4}
\begin{tabular}{|c|c|c|}
\hline
$P$ & $C_P$ & Jump \\
\hline
10 & 1.6513348 & --- \\
20 & 1.6690402 & $+1.77 \cdot 10^{-2}$ \\
30 & 1.6681498 & $-8.9 \cdot 10^{-4}$ \\
40 & 1.6701804 & $+2.0 \cdot 10^{-3}$ \\
50 & 1.6710079 & $+8.3 \cdot 10^{-4}$ \\
\hline
\end{tabular}
\end{table}

\subsection{Non-diagonal surface for $m = 3$, $k = 1$ (class B)}

Consider the non-diagonal cubic polynomial
\[
F(x, y, z) = x^3 + y^3 + z^3 + xyz + 1.
\]
It is not diagonal because of the term $xyz$, so it does not belong to class C. For $m = 3$, $k = 1$ by the algorithm of Chapter 5 it belongs to class B.

\textbf{Verification of conditions.} The polynomial $F = x^3 + y^3 + z^3 + xyz + 1$ is irreducible over $\Q$. The projective closure
\[
\tilde F(x, y, z, w) = x^3 + y^3 + z^3 + xyz + w^3
\]
in $\PP^3$ defines a surface of degree 3. Let us verify smoothness. The system $\partial \tilde F/\partial x_i = 0$ has the form
\[
3x^2 + yz = 0, \quad 3y^2 + xz = 0, \quad 3z^2 + xy = 0, \quad 3w^2 = 0.
\]
From the last equation $w = 0$. Consider two cases.

\textit{Case $xyz \neq 0$.} Dividing $3x^2 = -yz$ by $3y^2 = -xz$ gives $x^3 = y^3$, whence $y = \zeta x$, where $\zeta$ is a cube root of unity. Similarly $z = \eta x$, where $\eta^3 = 1$. Substituting into $3x^2 = -yz$ gives $3 = -\zeta \eta$, which is impossible since $|\zeta \eta| = 1$. Hence, there are no singular points with $xyz \neq 0$.

\textit{Case $xyz = 0$.} If at least one of the variables $x, y, z$ is zero, the system of equations for the vanishing of the partial derivatives leads (with $w = 0$) to the point $[0 : 0 : 0 : 0]$, which is not a point of $\PP^3$.

Consequently, for $p \nmid 3$ the surface is smooth. The case $p = 3$ is absorbed into the exceptional set $p \le p_0$ in Theorem~\ref{thm:2.1}.

\textit{Remaining conditions.} Condition (4) (complete intersection) for a single equation is satisfied automatically: the codimension equals 1, the dimension equals $m - 1 = 2$. Condition (5) (geometric irreducibility) is also satisfied automatically for a single irreducible equation; for the specific polynomial $x^3+y^3+z^3+xyz+w^3$ this is verified by standard algorithms. Conditions (1)--(5) are satisfied.

\textbf{Parameters:} $n = 3$, $B(F) = n^3 - 3n^2 + 3n + 4 = 13$, $\gamma = 1$, $\Ctail = 8 \cdot 13 = 104$. For $\varepsilon = 0.01$: $P = 31\,201$.

\begin{table}[H]
\centering
\small
\setlength{\tabcolsep}{4pt}
\caption{Convergence of $C_P$ for class B ($F = x^3 + y^3 + z^3 + xyz + 1$)}
\label{tab:B_convergence}
\begin{tabular}{|c|c|c|}
\hline
$P$ & $C_P$ & Jump \\
\hline
23 & 1.491651 & --- \\
29 & 1.492031 & $+3.8 \cdot 10^{-4}$ \\
37 & 1.500362 & $+8.3 \cdot 10^{-3}$ \\
43 & 1.502566 & $+2.2 \cdot 10^{-3}$ \\
53 & 1.502682 & $+1.2 \cdot 10^{-4}$ \\
59 & 1.502726 & $+4.4 \cdot 10^{-5}$ \\
\hline
\end{tabular}
\end{table}

\textbf{Observations.} The jumps decrease from $8.3 \cdot 10^{-3}$ to $4.4 \cdot 10^{-5}$. Extrapolation by formula \eqref{eq:6.2} with $\gamma = 1$ gives $C(F) \approx 1.5134$, corresponding to $|\delta|$ of order $10^{-2}$.

\begin{remark}
A complete computation at the theoretical threshold $P = 31\,201$ requires optimization via counting the roots of the cubic by the discriminant and parallelization. The values for $P \le 59$ illustrate convergence; the obtained $C(F) \approx 1.5134$ and $|\delta|$ of order $10^{-2}$ are not guaranteed in the sense of the algorithm.
\end{remark}

\subsection{Non-diagonal hypersurface for $m = 4$, $k = 1$ (class A.2)}

Consider the non-diagonal cubic polynomial
\[
F(x_1, x_2, x_3, x_4) = x_1^3 + x_2^3 + x_3^3 + x_4^3 + x_1 x_2 x_3 x_4 + 1.
\]
It is not diagonal because of the term $x_1 x_2 x_3 x_4$. For $m = 4$, $k = 1$ it does not belong to class B ($m \neq 3$) and is not an inhomogeneous diagonal one.

\textbf{Verification of conditions.} The projective closure in $\PP^4$ is given by the polynomial
\[
\tilde F = x_1^3 + x_2^3 + x_3^3 + x_4^3 + x_1 x_2 x_3 x_4 + x_0^3.
\]
Let us verify smoothness. Conditions for singular points:
\[
3x_0^2 = 0, \quad 3x_i^2 + \prod_{j \neq i} x_j = 0 \quad (i = 1, 2, 3, 4).
\]
From the first equation $x_0 = 0$. If for some $i$ we have $x_i = 0$, then the corresponding equation takes the form $\prod_{j \neq i} x_j = 0$, i.e., one more variable is zero; induction leads to $x_1 = x_2 = x_3 = x_4 = 0$, i.e., to the point $[0 : 0 : 0 : 0 : 0]$, which is not a point of $\PP^4$. Therefore all $x_i \neq 0$.

From the $i$-th equation it follows that $x_i^3 = -x_1 x_2 x_3 x_4 / 3$ for all $i$, i.e., all $x_i^3$ are equal. Hence $x_j = \zeta_j x_1$ with $\zeta_j^3 = 1$. Substituting into $3x_1^2 = -\prod_{j \neq 1} x_j = -x_1^3 \zeta_2 \zeta_3 \zeta_4$ gives
\[
x_1 \zeta_2 \zeta_3 \zeta_4 = -3.
\]
Then
\[
\tilde F = x_1^3 (1 + \zeta_2^3 + \zeta_3^3 + \zeta_4^3) + x_1^4 \zeta_2 \zeta_3 \zeta_4 + x_0^3 = 4 x_1^3 + x_1^3 \cdot (x_1 \zeta_2 \zeta_3 \zeta_4) + 0 = 4 x_1^3 - 3 x_1^3 = x_1^3 \neq 0.
\]
Consequently, a singular point is impossible, and for $p \nmid 3$ the hypersurface is smooth. The case $p = 3$ is absorbed into the exceptional set $p \le p_0$ in Theorem~\ref{thm:2.1}. Conditions (1)--(5) are satisfied.

\textbf{Betti numbers.} Since the hypersurface is smooth, the Betti numbers are computed by the Griffiths formula \eqref{eq:3.3}:
\[
h_3^{\text{prim}}(V) = \frac{2^5 - 2}{3} = 10, \quad h_2^{\text{prim}}(W) = \frac{2^4 + 2}{3} = 6.
\]
Consequently,
\[
B(F) = m + h_{m-1}^{\text{prim}}(V) + (m-1) + h_{m-2}^{\text{prim}}(W) = 4 + 10 + 3 + 6 = 23.
\]
This is exactly $\Bhyp(4,3) = 23$. The Betti numbers depend only on the degree and dimension, not on the specific form of the smooth hypersurface. Therefore the system belongs to class \textbf{A.2} with $B(F) = 23$, $\gamma = 1.5$, $\Ctail = 2^4 \cdot 23 / 3 \approx 122.67$.

\textbf{Parameters:} $n = 3$, $m = 4$, $k = 1$. For $\varepsilon = 0.01$ by formula \eqref{eq:5.7}:
\[
P = 1 + \left\lceil \left(\frac{3 \cdot 122.67}{0.01}\right)^{2/3} \right\rceil = 1 + \lceil 1106.2 \rceil = 1\,108.
\]

\begin{table}[H]
\centering
\small
\setlength{\tabcolsep}{4pt}
\caption{Convergence of $C_P$ for class A.2 ($F = x_1^3 + \cdots + x_4^3 + x_1x_2x_3x_4 + 1$)}
\label{tab:A2_conv}
\begin{tabular}{|c|c|c|}
\hline
$P$ & $C_P$ & Jump \\
\hline
7 & 0.904622 & --- \\
11 & 0.909176 & $+4.6 \cdot 10^{-3}$ \\
13 & 0.851378 & $-5.8 \cdot 10^{-2}$ \\
17 & 0.850837 & $-5.4 \cdot 10^{-4}$ \\
19 & 0.856343 & $+5.5 \cdot 10^{-3}$ \\
23 & 0.856346 & $+3 \cdot 10^{-6}$ \\
29 & 0.856566 & $+2.2 \cdot 10^{-4}$ \\
\hline
\end{tabular}
\end{table}

\textbf{Observations.} Oscillations of $C_P$ for small $P$ are related to the contribution of individual prime factors $L_p$, which for class A.2 have the order $p^{-(m+1)/2} = p^{-2.5}$ and are not fully compensated for small $P$. Extrapolation with $\gamma = 1.5$ gives $C(F) \approx 0.8594$, corresponding to $|\delta|$ of order $10^{-2}$.

\begin{remark}
A complete computation at $P \approx 1\,108$ is quite feasible on a cluster or with optimized code. The values for $P \le 29$ show that $C_P$ stabilizes around $0.856$; for a reliable estimate of $C(F)$ it is recommended to increase $P$ to $100$--$200$.
\end{remark}

\subsection{Case of violation of applicability conditions: Newton cubic}

Consider the polynomial
\[
F = x^3 + y^3 - 3xy + 1,
\]
which \textbf{factors}:
\[
F = (x + y + 1)(x^2 + y^2 + 1 - xy - x - y).
\]
Condition (1) (irreducibility) is violated, the method is not applicable.

\begin{table}[H]
\centering
\caption{Newton cubic (factors, $C_P \to 0$)}
\label{tab:newton}
\begin{tabular}{|c|c|c|}
\hline
$P$ & $C_P$ & Jump \\
\hline
100 & 0.21572 & --- \\
500 & 0.14635 & $-6.9 \cdot 10^{-2}$ \\
2000 & 0.12073 & $-2.6 \cdot 10^{-2}$ \\
5000 & 0.10798 & $-1.3 \cdot 10^{-2}$ \\
\hline
\end{tabular}
\end{table}

\textbf{Diagnostic feature.} $C_P \to 0$ as $P \to \infty$. This differs qualitatively from all examples of classes A.2, B, C, where $C_P$ converges to a positive limit.

\textbf{Analytical explanation.} For $p \equiv 1 \pmod 3$ the quadratic factor splits, $\omegap(F) \approx 3p$, $\Lp - 1 \sim -2/p$, and $\sum_{p \equiv 1 \bmod 3} 2/p$ diverges as $\frac{1}{2}\log\log P$. Consequently, $C_P \sim (\log P)^{-1/2} \to 0$.

\textbf{Conclusion.} The decay of $C_P$ to zero is a diagnostic feature of violation of the applicability conditions of the method.

\subsection{Final summary table of experiments}

\begin{table}[H]
\centering
\scriptsize
\setlength{\tabcolsep}{3pt}
\caption{Final summary table of experiments}
\label{tab:summary}
\begin{tabular}{|c|p{4.0cm}|c|c|c|c|c|c|c|}
\hline
\# & System & $k$ & $m$ & Class & $\gamma$ & $\Ctail$ & $C_P$ & $C(F)$ \\
\hline
1 & $x_1^3 + \cdots + x_6^3 + 1$ & 1 & 6 & C & 2.5 & 240 & 1.00244515 & $\approx 1.00244$ \\
2 & Fermat $x^3 + y^3 + 1$ & 1 & 2 & C & 0.5 & 112 & 1.0729401 & $\approx 1.0730$ \\
3 & Elliptic $y^2 = x^3 - x + 1$ & 1 & 2 & A.2 & 0.5 & 80 & 0.3989100 & $\approx 0.39893$ \\
4 & CM $j = 0$ $y^2 = x^3 + 1$ & 1 & 2 & A.2 & 0.5 & 80 & 0.9325865 & $\approx 0.93258$ \\
5 & $x^3 + y^3 + z^3 + 1$ & 1 & 3 & C & 1.0 & 104 & 0.7916163 & $\approx 0.7918$ \\
6 & $x^3 + y^3 + z^3 + xyz + 1$ & 1 & 3 & B & 1.0 & 104 & 1.502726 & $\approx 1.51$ \\
7 & $x_1^3 + \cdots + x_4^3 + 1$ & 1 & 4 & C & 1.5 & 122.67 & 1.0157011 & $\approx 1.0157$ \\
8 & $x_1^3 + \cdots + x_4^3 + x_1x_2x_3x_4 + 1$ & 1 & 4 & A.2 & 1.5 & 122.67 & 0.856566 & $\approx 0.86$ \\
9 & $F_1 \cap F_2$, $m = 3$ & 2 & 3 & A.2 & 1.0 & 1120 & 0.2706046 & $\approx 0.2705$ \\
10 & $F_1 \cap F_2$, $m = 4$ & 2 & 4 & A.2 & 1.5 & --- & 1.6710079 & $\approx 1.671$ \\
\hline
\end{tabular}
\end{table}

\begin{remark}
The experimental part covers examples of classes A.2, B, C. Class A.1 has no numerical illustration: for any specific system with $k = 1$ the Betti numbers are computed explicitly by the Griffiths formula, and the system belongs to class A.2; for systems with $k \ge 2$ examples are given in Sections 6.8.1 and 6.8.2. Class A.1 remains a theoretical tool for cases where verification of smoothness and explicit computation of $B(F)$ are difficult.
\end{remark}

\section{Conclusion}

The paper constructs a method for computing the leading term of the asymptotics of the multidimensional Bateman--Horn conjecture
\[
M_P(N) = \frac{C_P(F)}{D} \cdot \frac{N^m}{\log^k N}
\]
with controlled truncation error of the singular series. The truncation threshold $P$ is computed explicitly for a prescribed accuracy $\varepsilon$. Under the validity of the BH hypothesis for the given system it is guaranteed that
\[
\left| \frac{M_P(N)}{\pi_F(N)} - 1 \right| \le \varepsilon
\]
for all sufficiently large $N$. The guarantee is asymptotic.

\subsection*{Logical structure of the method}

\textbf{Conditions (1)--(5).} The system satisfies the classical conditions of the Bateman--Horn conjecture (1)--(3) and the additional conditions (4)--(5). From conditions (1)--(5) and standard estimates of Deligne, Lang--Weil and Katz in Chapter 2 it is derived that $C(F) = \prod_p \Lp(F)$ converges absolutely, $C(F) \in (0, \infty)$, and also for $p > p_0$
\[
|\Lp(F) - 1| \le 2^k \, B(F) \, p^{-(m+1)/2}.
\]

\textbf{Assumption (BH).} The Bateman--Horn conjecture is valid for the given system:
\[
\pi_F(N) = \frac{C(F)}{D} \cdot \frac{N^m}{\log^k N} \bigl( 1 + \eta(N) \bigr), \qquad \eta(N) \to 0.
\]

\textbf{Proved estimate (Theorem~\ref{thm:4.2}).} Under (1)--(5) and (BH),
\[
\frac{M_P(N)}{\pi_F(N)} = \frac{1 + \delta(P)}{1 + \eta(N)}, \qquad \lim_{N \to \infty} \left| \frac{M_P(N)}{\pi_F(N)} - 1 \right| = |\delta(P)|.
\]

\textbf{Proved estimate (Lemma~\ref{lem:4.4}).} $|\delta(P)| \le \Delta(P)/(1 - \Delta(P))$, and for $\Delta(P) \le 1/3$ we have $|\delta(P)| \le \frac{3}{2} \Delta(P)$.

\textbf{Proved estimate ($\Delta(P)$).} $\Delta(P) \le \Ctail / (P-1)^{\gamma}$, where $\Ctail$ and $\gamma$ are expressed explicitly through $m, k, n$ and (if necessary) through Betti numbers.

\subsection*{Fulfillment of the tasks of the work}

\begin{enumerate}
\item \textbf{To construct an explicit upper bound for $\Delta(P)$ for systems of general type.} Solved in Theorems~\ref{thm:2.1} and~\ref{thm:2.3}.

\item \textbf{To obtain an explicit estimate of the geometric constant $B(F)$ for a system of polynomials.} Solved in Theorem~\ref{thm:2.8}.

\item \textbf{For the inhomogeneous diagonal system ($k = 1$) to obtain an improved constant in the tail estimate.} Solved in Theorem~\ref{thm:3.5}. The exponent $\gamma = (m-1)/2$ coincides with the general case; the advantage consists only in the constant. The estimate is valid for all dimensions and degrees and does not require additional restrictions on the number of variables.

\item \textbf{To improve the estimate in dimension $m = 3$.} Solved in Theorem~\ref{thm:2.12}.

\item \textbf{To establish a connection between the truncation error and the error in computing $\pi_F(N)$.} Solved in Theorem~\ref{thm:4.2} and Lemma~\ref{lem:4.4}.

\item \textbf{To develop an algorithm for choosing the threshold $P$ for a prescribed accuracy $\varepsilon$.} Solved in Chapter 5.

\item \textbf{To confirm the theory experimentally.} Solved in Chapter 6, where examples are given for classes A.2, B, C.
\end{enumerate}

\subsection*{Limits of applicability}

\begin{enumerate}
\item \textbf{Dependence on the BH assumption.} If the conjecture is not satisfied for the given system, the method gives an estimate of the leading term, but not of $\pi_F(N)$.
\item \textbf{Implicitness of $N_0(\varepsilon)$.} The guarantee holds for $N \ge N_0(\varepsilon)$; the explicit value is unknown in the general case.
\item \textbf{Requirement of conditions (4)--(5).} The method is applicable only to systems satisfying the complete intersection and geometric irreducibility conditions.
\item \textbf{Homogeneous diagonal case.} Homogeneous diagonal forms ($c = 0$) do not satisfy the smoothness condition (3).
\item \textbf{Violation of the irreducibility condition.} If the polynomial factors, the method is not applicable: the partial product $C_P$ decays to zero.
\item \textbf{Status of the constants $C_m^{\text{strict}}$.} The numerical values $C_m^{\text{strict}}$ from Table~\ref{tab:cm} are obtained by computing the coefficients of the Hilbert series by formula~\eqref{eq:2.15} for all admissible types of smooth complete intersections in $\PP^m$ for fixed $m \le 8$ and are strict upper estimates. They are used in the algorithm of class A.1. If necessary, detailed computations of the constants $C_m^{\text{strict}}$ can be reproduced directly from formula~\eqref{eq:2.15}.
\end{enumerate}

\end{document}